\documentclass[12pt,a4paper,leqno]{article}

\usepackage{amssymb,exscale}
\usepackage[centertags]{amsmath}
\usepackage{amsthm}

\numberwithin{equation}{section}

\swapnumbers
\theoremstyle{definition}
\newtheorem{theorem}[equation]{Theorem}
\newtheorem{lemma}[equation]{Lemma}
\newtheorem{proposition}[equation]{Proposition}
\newtheorem{definition}[equation]{Definition}
\newtheorem{example}[equation]{Example}

\newtheorem{myItem}[equation]{}

\renewcommand{\phi}{\varphi}

\newcommand{\ti}{\tilde}

\renewcommand{\(}{\bigl(}
\renewcommand{\)}{\bigr)\vphantom{)}}

\newcommand{\Ind}{\operatorname{Ind}}
\newcommand{\mes}{\operatorname{L}}

\renewcommand{\Im}{\operatorname{Im}}
\newcommand{\const}{\operatorname{const}}

\newcommand{\al}{\alpha}
\newcommand{\eps}{\varepsilon}
\newcommand{\ga}{\gamma}
\newcommand{\si}{\sigma}
\newcommand{\cS}{{\mathfrak S}}
\newcommand{\F}{\mathcal F}
\newcommand{\G}{\mathcal G}
\renewcommand{\P}{\mathcal P}
\newcommand{\la}{\lambda}
\newcommand{\R}{\mathbb R}
\newcommand{\C}{\mathbb C}
\newcommand{\Z}{\mathbb Z}

\newcommand{\equi}{\;\;\Longleftrightarrow\;\;}

\def\emailwww#1#2{\par\qquad {\tt #1}\par\qquad {\tt #2}\medskip}

\newenvironment{myitemize}{\begin{list}{$\bullet$}
{\setlength{\topsep}{1mm}
\setlength{\partopsep}{0mm}
\setlength{\itemsep}{1mm}
\setlength{\parsep}{0mm}
\setlength{\parskip}{0mm}}}
{\end{list}}

\begin{document}

\title{Spectral densities describing off-white noises}
\author{Boris Tsirelson}
\date{}

\maketitle

\stepcounter{footnote}
\footnotetext{Supported in part by the Israel Science Foundation grant
592/99/2.}

\begin{abstract}
For the white noise, the spectral density is constant, and the past
(restriction to $ (-\infty,0) $) is independent from the future
(restriction to $ (0,+\infty) $). If the spectral density is not too
far from being constant, then dependence between the past and the
future can be eliminated by an equivalent measure change; that is
called an off-white noise. I derive from well-known results a
necessary and sufficient condition for a spectral density to describe
an off-white noise.
\end{abstract}

\section*{Introduction}

`Past and future' is a well-known topic in the theory of stationary
Gaussian random processes. The restriction of a process $ X(t) $ to $
t \in (-\infty,0) $ is the past; the future is its restriction to $
(0,+\infty) $ or, more generally, $ (s,+\infty) $. Typically one shows
that the past and the future are nearly independent if the separation
$ s $ is large enough, under appropriate conditions on the spectral
density of the process. In contrast, the present work deals with the
case $ s = 0 $ (no separation). For a continuous process, of course, $
X(0) $ belongs both to the past and to the future, making them heavily
dependent. However, for the white noise they are independent
anyway. An off-white noise\footnote{%
 The term is suggested by William Arveson; I called it a `slightly
 coloured noise'.}
introduced in \cite{Ts} (motivated by the theory of continuous tensor
products of Hilbert spaces) generalizes the white noise. It is defined
as a stationary Gaussian generalized random process such that the
joint distribution of the past and the future is equivalent (that is,
mutually absolutely continuous) to the product of their marginal
distributions. The present work derives from well-known results
about `past and future' a necessary and sufficient condition for a
spectral density to describe an off-white noise. I feel that it is
basically a folklore worth to be written.

I am indebted to Ildar Ibragimov, Alexander Kheifets, Mikhail Sodin
and Sergei Treil; due to their patient efforts I am now less afraid of
the theory of analytic functions.

\vspace{-3mm}
\hfill\hrulefill\hfill\mbox{}
\vspace{1mm}

In the Hilbert space $ H $ of all linear measurable functionals over a
Gaussian random process, the past $ \P $ and the future $ \F $ are
linear subspaces. For the white noise, $ H = \P \oplus \F $, the
subspaces being orthogonal. For an off-white noise the corresponding
relation is
\[
H = \P \oplus \F \quad \text{in the FHS sense,}
\]
as defined in \cite{Ts}; it means that the orthogonal projection from
$ \F $ to $ \P $ is a Hilbert-Schmidt operator, and $ H = \P \oplus \F
$ in the topological sense. The latter means that $ \P + \F $ is dense
in $ H $ and $ \P,\F $ are at positive angle (that is, the projection
is of norm less than $ 1 $).

\section{Analytic functions inside and outside the circle}

\begin{definition}\label{def1}
(a) A \emph{past-and-future structure} (or `PaF structure') consists
of:
\begin{myitemize}
\item a separable Hilbert space $ H $;
\item a two-sided sequence $ (\P_n)_{n\in\Z} $ of (closed linear)
 subspaces $ \P_n \subset H $, increasing (in the sense that $ \P_n
 \subset \P_{n+1} $ for all $ n $) and such that the union of all $
 \P_n $ is dense in $ H $;
\item a two-sided sequence $ (\F_n)_{n\in\Z} $ of subspaces $ \F_n
 \subset H $, decreasing and such that the union of all $ \F_n $ is
 dense in $ H $;
\item a unitary operator $ T : H \to H $ such that $ T \P_n = \P_{n+1}
 $ and $ T \F_n = \F_{n+1} $ for all $ n $.
\end{myitemize}

(b) Two PaF structures $ ( H, (\P_n), (\F_n), T ) $ and $ ( H',
(\P'_n), (\F'_n), T' ) $ are \emph{isomorphic,} if there exists an
invertible linear isometry $ U : H \to H' $ such that $ U \P_n =
\P'_n $ and $ U \F_n = \F'_n $ for all $ n $, and $ U T = T' U
$.

(c) A \emph{PaF geometry} is a PaF structure treated up to
isomorphism.

(d) For any PaF geometry $ \G = ( H, (\P_n), (\F_n), T ) $ and any $ k
\in \Z $ define $ \G+k $ (the shifted PaF geometry) as $ ( H,
(\P_{n+k}), (\F_n), T ) $ (or equivalently $ ( H, (\P_n),
(\F_{n-k}), T ) $). Also define the \emph{time-reversed PaF geometry}
as $ \overline\G = ( H, (\F_{-n}), (\P_{-n}), T^{-1} ) $.
\end{definition}

Let $ \mu $ be a (positive) $ \sigma $-finite Borel measure on the
unit circle $ \{ z \in \C : |z| = 1 \} $. The set of all polynomials $
P $ such that $ \int |P|^2 \, d\mu < \infty $ is an ideal in the
commutative ring of polynomials. If the ideal contains not only $ 0 $,
then it is generated by a single polynomial $ P_\mu $ (not identically
$ 0 $), since every ideal in that ring is principal. It is easy to see
that (up to a coefficient that may be ignored), $ P_\mu (z) = (z-z_1)
\dots (z-z_m) $ for some $ z_1, \dots, z_m $ on the circle. If $ \mu $
is finite then $ m = 0 $ and $ P_\mu (z) = 1 $.

\begin{definition}\label{def2}
(a) A nonatomic $ \sigma $-finite Borel measure $ \mu $ on the circle
will be called \emph{moderate,} if the ideal of polynomials $ P $
satisfying $ \int |P|^2 \, d\mu < \infty $ is different from $ \{ 0
\} $.

(b) \emph{Poles} of a moderate measure $ \mu $ are roots of the
polynomial $ P_\mu $ that generates the ideal. \emph{Multiplicity} of
a pole of $ \mu $ is its multiplicity as a root of $ P_\mu $.
\end{definition}

The set of all moderate
measures is a linear space, closed under multiplication by functions
of the form $ f/|P|^2 $ where $ f $ is a bounded nonnegative Borel
function on the circle and $ P $ is a polynomial (not identically $ 0
$).

Conjugation $ z \mapsto \overline z $ maps the circle onto itself, and
sends each measure $ \mu $ on the circle to another measure, denote it
$ \ti\mu $. Clearly, $ \mu $ is moderate if and only if $ \ti\mu $ is
moderate, and if they are, then $ P_\mu (z) = \overline{ P_{\ti\mu}
(\overline z) } $ for all $ z $. Also, each function $ f \in L_2 (\mu)
$ corresponds to another function $ \ti f \in L_2 (\ti\mu) $ such that
$ f(z) = \ti f (\overline z) $. However, if $ f $ is a polynomial $ P
$ restricted to the circle then $ \ti f $ is rather the rational
function $ z \mapsto P(1/z) $ restricted to the circle. In particular,
if $ P(z) = z - z_1 $ where $ |z_1| = 1 $ then $ P (1/z) = -z_1
(z-\overline z_1) / z = -z_1 \overline{ P (\overline z) } / z
$. Accordingly, if $ P(z) = (z-z_1) \dots (z-z_m) $ for some $ z_1,
\dots, z_m $ on the circle then $ P (1/z) = \const \cdot z^{-m}
\overline{ P (\overline z) } $. Therefore
\[
P_{\ti\mu} \Big( \frac1z \Big) = \const \cdot z^{-m} P_\mu (z) \quad
\text{for all } z \in \C \setminus \{0\} \, ; \quad \text{here } m =
\deg P_\mu \, ;
\]
and $ \ti P_{\ti\mu} (z) = \const \cdot z^{-m} P_\mu (z) $ for $ |z| =
1 $. Functions $ P_\mu $ and $ \ti P_{\ti\mu} $ have the same zeros
(on the circle); however, $ P_\mu $ has a pole (of multiplicity $ m $)
at $ \infty $, while $ \ti P_{\ti\mu} $, or rather its analytic
continuation $ \const \cdot z^{-m} P_\mu (z) $, has a pole (of
multiplicity $ m $) at $ 0 $.

Every moderate measure $ \mu $ determines a PaF structure $
( H, (\P_n), (\F_n), T ) $ as follows:
\begin{myitemize}
\item $ H = L_2 (\mu) $;
\item $ \F_n $ is spanned by functions $ z \mapsto z^k P_\mu (z) $
 for $ k \in \Z $, $ k \ge n \, $;\footnote{%
  Here $ P_\mu $ may be replaced with any polynomial $ P $ (not
  identically $ 0 $) satisfying $ \int |P|^2 \, d\mu < \infty $;
  redundant roots of $ P $ do not influence the \emph{closed}
  subspace.}
\item $ \P_n $ is spanned by functions $ z \mapsto z^k \ti P_{\ti\mu}
(z) $ for $ k \in \Z $, $ k \le n $;\item $ (Tf) (z) = z f(z) $ for $
f \in L_2(\mu) $.
\end{myitemize}

In other words, $ \P_n $ is spanned by functions $ z \mapsto z^k P_\mu
(z) $ for $ k \in \Z $, $ k \le n - \deg P_\mu $.

Treating the PaF structure up to isomorphism, we get a PaF geometry;
denote it by
\[
\G_\mu \, .
\]
The time-reversed PaF geometry (as defined by \ref{def1}(d))
corresponds to $ \ti\mu $:
\[
\overline{ \G_\mu } = \G_{\ti\mu} \, ;
\]
an isomorphism is $ L_2 (\mu) \ni f \mapsto \ti f \in L_2 (\ti\mu) $.
If $ \mu $ is symmetric (that is, $ \mu = \ti\mu $) then $ \G_\mu $ is
time-symmetric (that is, $ \G_\mu = \overline{\G_\mu} $).

\begin{proposition}\label{prop1}
Let $ \mu, \mu' $ be moderate measures, $ z_0 $ a point on
the circle, and $ \mu' (dz) = |z-z_0|^2 \mu (dz) $. Then\footnote{%
 $ \G_\mu + 1 $ is the shifted PaF geometry, recall \ref{def1}(d).}
\[
\G_{\mu'} = \G_\mu + 1 \, .
\]
\end{proposition}

\begin{proof}
We have $ \G_\mu = ( H, (\P_n), (\F_n), T ) $, $ \G_{\mu'} = ( H',
(\P'_n), (\F'_n), T' ) $.
Multiplication by $ 1 / (z-z_0) $ is an invertible linear isometry $
L_2 (\mu) \to L_2 (\mu') $, that is, $ H \to H' $; it intertwines $ T
$ with $ T' $. We'll prove that it sends $ \F_0 $ to $ \F'_0 $ and $
\P_1 $ to $ \P'_0 $.

First, $ \frac1{z-z_0} \F_0 \supset \F'_0 $ for a trivial reason: for
all $ k \ge 0 $ the function $ z \mapsto (z-z_0) z^k P_{\mu'} (z) $
belongs to $ \F_0 $, since $ \int | (z-z_0) P_{\mu'} (z) |^2 \,
\mu(dz) = \int | P_{\mu'} (z) |^2 \, \mu'(dz) < \infty $.\footnote{%
 It may happen that $ z_0 $ is a pole of $ \mu $, then $ P_\mu (z) =
 (z-z_0) P_{\mu'} (z) $; otherwise $ P_\mu = P_{\mu'} $. In any case
 $ z \mapsto (z-z_0) P_{\mu'} (z) $ belongs to the ideal generated by $
 P_\mu $.}

In order to prove that $ \frac1{z-z_0} \F_0 \subset \F'_0 $ take
polynomials $ P_n $ such that $ P_n (z) \to \frac1{z-z_0} $ for $ n
\to \infty $, and $ |P_n(z)| \le 2 \big| \frac1{z-z_0} \big| $,
whenever $ |z| = 1 $; say, we may take
\[
P_n (z) = - (1-\eps_n) \overline z_0 \frac{ 1 - (1-\eps_n)^n \overline
z_0^n z^n }{ 1 - (1-\eps_n) \overline z_0 z }
\]
choosing $ \eps_n \to 0+ $ such that $ n \eps_n \to \infty $. We have
(for every $ k \ge 0 $) $ P_n(z) z^k P_\mu (z) \to \frac1{z-z_0} z^k
P_\mu (z) $ pointwise, and $ | P_n(z) z^k P_\mu (z) | \le
\frac2{|z-z_0|} | P_\mu(z) | $. The majorant belongs to $ L_2 (\mu')
$; polynomials $ z \mapsto P_n (z) z^k P_\mu (z) $ belong to $ L_2
(\mu) $, therefore to $ L_2 (\mu') $, and to $ \F'_0 $. So, $
\frac1{z-z_0} \F_0 = \F'_0 $.

Now we apply the equality $ \frac1{z-z_0} \F_0 = \F'_0 $ to measures $
\ti\mu $, $ \ti\mu' $ (symmetric to $ \mu $, $ \mu' $); these are
related by $ \ti\mu'(dz) = | z - \overline z_0 |^2 \ti\mu (dz) $;
thus, $ \frac1{z-\overline z_0} \F_0 (\ti\mu) = \F_0 (\ti\mu') $. The
isomorphism $ f \mapsto \ti f $ between $ \G_\mu $ and $ \G_{\ti\mu} $
(as well as $ \G_{\mu'} $ and $ \G_{\ti\mu'} $) transforms $ \F_0
(\ti\mu) $ to $ \P_0 (\mu) $, $ \F_0 (\ti\mu') $ to $ \P_0 (\mu') $,
and the function $ z \mapsto \frac1{z-\overline z_0} $ into the
function $ z \mapsto \frac1{\overline z - \overline z_0 } $. So,
\[
\frac1{\overline z - \overline z_0 } \P_0 (\mu) = \P_0 (\mu') \, .
\]
However, $ \frac1{\overline z - \overline z_0 } = \const \cdot z \cdot
\frac1{z-z_0} $ for $ |z|=1 $ (namely, $ \const = -\overline z_0 $);
therefore
\[
\frac1{\overline z - \overline z_0 } \P_0 (\mu) = \frac1{z-z_0} z \P_0
(\mu) = \frac1{z-z_0} \P_1 (\mu) \, .
\]
So, $ \frac1{z-z_0} \P_1 (\mu) = \P_0 (\mu') $, that is, $
\frac1{z-z_0} \P_1 = \P_0' $.
\end{proof}

Given a PaF geometry $ \G = ( H, (\P_n), (\F_n), T ) $, we may ask,
whether or not two spaces $ \P_n, \F_{n+k+1} $ are at positive
angle.\footnote{%
 Alternatively we could ask whether or not they are orthogonal, have trivial
 intersection, etc. Every such property leads to its `index' satisfying
 \eqref{eqind}.}
It depends on $ k $, not $ n $. If it holds for $ k $ then it surely
holds for $ k+1 $. We define the index, $ \Ind(\G) $, as the least $ k
\in \Z $ possessing the property.\footnote{%
 If all $ k $ possess the property then $ \Ind(\G) = -\infty $; if
 no one does then $ \Ind(\G) = +\infty $.}
Evidently,
\begin{equation}\label{eqind}
\Ind ( \G + k ) = \Ind (\G) + k \, .
\end{equation}
Combined with Proposition \ref{prop1} it means that $ \Ind (G_{\mu'})
= \Ind (G_\mu) + 1 $ whenever $ \mu' (dz) = |z-z_0|^2 \mu (dz) $.

Assume for a while that $ \mu $ is finite. We have $ \Ind (\G_\mu) \ge
0 $, since constant functions belong both to $ \P_0 $ and to $ \F_0 $.
It is well-known (see \cite[Sect.~9]{HS} or \cite[Th.~4 in
Sect.~V.2]{IR}) that $ \Ind (\G_\mu) = N $ if and only if $ d\mu =
|P|^2 \, d\nu $ for some polynomial $ P $ of degree $ N $ with all
roots on the circle, and some \emph{finite} measure $ \nu $ such that
$ \Ind (\G_\nu) = 0 $. Finiteness of $ \nu $ ensures that $
\P_0^{(\nu)} + \F_1^{(\nu)} $ is dense in $ H^{(\nu)} $. Thus $
H^{(\nu)} = \P_0^{(\nu)} \oplus \F_1^{(\nu)} $ in the topological
sense. Taking into account that $ \G_\mu = \G_\nu + N $ we see that
the two following conditions are equivalent for every finite measure $
\mu $:
\begin{subequations}
\begin{gather}
\Ind (\G_\mu) = N \, ; \label{1.5a} \\
H = \P_0 \oplus \F_{N+1} \text{ in the topological sense.}
 \label{1.5b}
\end{gather}
\end{subequations}
Therefore (due to \ref{prop1}) these conditions are equivalent for
every \emph{moderate} measure $ \mu $.

In order to get $ H = \P_0 \oplus \F_{N+1} $ in the FHS sense, (one
of) the following two equivalent conditions must be added:
\begin{subequations}
\begin{myItem}
The orthogonal projection from $ \F_{N+1} $ to $ \P_0 $ is a
 Hilbert-Schmidt operator.\label{1.6a}
\end{myItem}
\begin{myItem}The product $ P_0 F_{N+1} P_0 $ is a trace-class operator;
 here $ P_0 $ and $ F_{N+1} $ are orthogonal projections (from $ H $)
 to $ \P_0 $ and $ \F_1 $ respectively.\label{1.6b}
\end{myItem}
\end{subequations}

Recall that a real-valued function $ \phi $ on the circle belongs to
Sobolev space $ W_2^{1/2} $ if and only if it satisfies the following
two equivalent conditions:
\begin{subequations}
\begin{equation}\label{eqSobolev}
\iint \frac{ | \phi(z_1) - \phi(z_2) |^2 }{ | z_1 - z_2 |^2 } \mes
(dz_1) \mes (dz_2) < \infty \, ,
\end{equation}
where $ \mes $ stands for Lebesgue measure;
\begin{equation}\label{eqSobolevF}
\sum_{-\infty}^{+\infty} |n| |\hat\phi_n|^2 < \infty \, ,
\end{equation}
where $ \hat\phi_n $ are Fourier coefficients of $ \phi $.
\end{subequations}

A well-known deep result of Ibragimov and Solev (see
\cite[Sect.~IV.4]{IR}, see also \cite[Sect.~7]{Pe}) states that a 
finite measure $ \mu $ satisfies \emph{both} \eqref{1.5a} \emph{and}
\eqref{1.6b} if and only if $ \mu $ has a density $ w $ (w.r.t.\
Lebesgue measure) of the form $ w = |P|^2 \exp \phi $ where $ \phi \in
W_2^{1/2} $ and $ P $ is a polynomial of degree $ N $ with all zeros
on the circle.

Combining the deep result with Proposition \ref{prop1} we generalize
the former from finite to moderate measures as follows.

\begin{proposition}\label{prop2}
For every moderate measure $ \mu $ on the circle and integer $ N $,
the following two conditions are equivalent.

(a) $ H = \P_0 \oplus \F_{N+1} $ in the FHS sense.

(b) $ \mu $ has a density of the form $\displaystyle
\frac{d\mu}{d\mes} (z) = \left| \frac{ (z-z_1)
\dots (z-z_l) }{ (z-z'_1) \dots (z-z'_m) } \right|^2 \exp \phi(z) $
for some $ l,m \in \{0,1,2,\dots\} $ such that $ l-m = N $, some points
$ z_1, \dots, z_l $, $ z'_1, \dots, z'_m $ on the circle, and some
function $ \phi \in W_2^{1/2} $.
\end{proposition}

\begin{proof}
A moderate measure $ \mu $ is related to a finite measure $ \nu $ by $
d\nu = |P_\mu|^2 \, d\mu $; thus $ \G_\nu = \G_\mu + m $ where $ m =
\deg P_\mu $. Condition (a) for $ \mu $ is equivalent to the condition
$ H^{(\nu)} = \P_0^{(\nu)} \oplus \F_{N+m+1}^{(\nu)} $ for $ \nu
$. The latter holds if
and only if $ \nu $ has a density $ w $ of the form $ w = |P|^2
\exp\phi $, where $ \phi \in W_2^{1/2} $ and $ \deg P = N + m $. It
means that $ \mu $ has the density $ \frac{|P|^2}{|P_\mu|^2} \exp\phi
$; note that $ \deg P - \deg P_\mu = N $.
\end{proof}

The following remarks will not be used.

If $ \mu $ satisfies Condition \ref{prop2}(b) then $ l,m $ and $
z_1,\dots,z_l $, $ z'_1,\dots,z'_m $ are uniquely determined by $ \mu
$ (provided that $ z_i \ne z'_j $ for all $ i,j $, of course).

A proposition similar to \ref{prop2} holds for ``$ H = \P_0 \oplus
\F_{N+1} $ in the
topological sense''; here the condition ``$ \phi \in W_2^{1/2} $'' is
replaced with the Helson-Szeg\"o condition: $ \phi = \tilde\psi + \chi
$ with $ \| \psi \|_\infty < \frac{\pi}2 $ and $ \| \chi \|_\infty <
\infty $, where $ \psi, \chi $ belong to $ L_\infty $ on the circle,
and $ \tilde\psi $ is the conjugate function to $ \psi $. (Or
alternatively, Muckenhoupt's condition $ (A_2) $ may be used.)

\section{Generalized random processes in continuous time}

Consider a Gaussian measure $ \ga $ in the space of (tempered,
Schwartz; real-valued) distributions (generalized functions) over $ \R
$; assume that $ \ga $ is invariant under shifts of $ \R $. Such
measures are probability distributions\footnote{%
 Sorry, `a distribution in the space of distributions' may be
 confusing. A `probability distribution' is just a probability measure
 (intended to describe a random element of the corresponding space). In
 contrast, a generalized function, called also `distribution', is a
 more singular (than a measure) object over $ \R $, generally not
 positive; for example, a derivative $ \delta^{(n)} $ of Dirac's
 delta-function.}
of stationary Gaussian generalized random
processes \cite{Ito}. The space of tempered distributions is dual to
the space of rapidly
decreasing infinitely differentiable functions $ \phi $ on $ \R $. 
Such $ \phi $ gives a linear functional on the space of distributions;
w.r.t.\ $ \ga $ it gives a normally distributed random variable, whose
variance is a quadratic form of $ \phi $ and may be written as $ \int
|\hat\phi|^2 \, d\nu $ where $ \hat\phi $ is Fourier transform of $
\phi $, and $ \nu $ is so-called spectral measure (of $ \ga $). It is a
positive $ \si $-finite Borel measure on $ \R $, symmetric (that is,
invariant under the map $ \la \to -\la $) and such that
\begin{equation}\label{2.1}
\int_{-\infty}^{+\infty} \frac1{(1+\la^2)^m} \, \nu(d\la) < \infty
\end{equation}
for $ m $ large enough, see \cite[Th.~3.3, 3.4]{Ito}.
Let $ m \in \{ 0,1,2,\dots \} $ be the least number satisfying
\eqref{2.1}. If $ \nu $ is finite then $ m = 0 $.

Consider the Sobolev space $ W_1^m (\R) $ of all functions $ \phi :
\R \to \R $ such that $ \phi, \phi', \dots, \phi^{(m)} \in L_1 (\R)
$. If $ \phi \in W_1^m (\R) $ then functions $ \la \mapsto \hat\phi
(\la) $, $ \la \mapsto \la \hat\phi (\la) $, \dots, $ \la \mapsto
\la^m \hat\phi (\la) $ belong to the space $ C_0 (\R) $ of all
(bounded) continuous functions on $ \R $ vanishing at $ \infty $, which
means that the function $ \lambda \mapsto (1+\la^2)^{m/2} \hat\phi
(\la) $ belongs to $ C_0 (\R) $. Taking into account \eqref{2.1} we have
\[
\forall \phi \in W_1^m (\R) \quad \hat\phi \in L_2 (\nu) \, .
\]
Thus, the quadratic form $ \phi \mapsto \int |\hat\phi|^2 \,
d\nu $ extends naturally from the space of rapidly increasing
infinitely differentiable functions to $ W_1^m (\R) $. Of course,
the former space is dense in the latter.\footnote{%
 Note also that $ C_0 (\R) $ could be replaced with $ C (\R) $ (all
 bounded continuous functions on $ \R $); accordingly, $ \phi^{(m)} $
 could be a finite measure rather than a function of $ L_1 (\R) $,
 which will be used in the proof of Lemma \ref{2.2}.}

Introduce two subspaces $ \P_0 (\nu), \F_0 (\nu) \subset L_2 (\nu) $;
namely, $ \P_0 (\nu) $ is spanned by functions $ \hat\phi $ where $
\phi \in W_1^m (\R) $, $ \phi(t) = 0 $ for $ t \in [0,\infty) $;
the same for $ \F_0 (\nu) $, but $ \phi(t) = 0 $ for $ t \in
(-\infty,0] $. The map $ \la \mapsto -\la $ sends $ \nu $ to itself,
and $ \P_0 (\nu) $ to $ \F_0 (\nu) $. That is, $ f \in \F_0 (\nu) $ if
and only if $ \ti f \in \P_0 (\nu) $; here $ \ti f (-\la) = f (\la) $.

We use the well-known conformal map $\displaystyle z = \frac{ z' - i
}{ z' + i } $, $\displaystyle z' = -i \frac{z+1}{z-1} $ of the real
line $ \Im z' = 0 $ to the unit circle $ |z| = 1 $; it also maps
half-planes $ \Im z' > 0 $, $ \Im z' < 0 $ onto the disk $ |z| < 1 $
and the region $ |z| > 1 $ respectively. Denote by $ \mu $ the image
of $ \nu $ under the map $ \la \mapsto \frac{\la-i}{\la+i} $; the $
\si $-finite measure $ \mu $ on the circle is symmetric (that is, $
\mu = \ti\mu $) and satisfies
\[
\int |1-z|^{2m} \, \mu(dz) < \infty \, ,
\]
which is the same as \eqref{2.1}, since $ \big| 1 -
\frac{\la-i}{\la+i} \big|^2 = \frac4{\la^2+1} $. In terms of
Definition \ref{def2}, $ \mu $ is a moderate measure; it has a pole of
multiplicity $ m $ at $ 1 $ (or it is a finite measure, and $ m=0 $);
$ P_\mu (z) = (z-1)^m $; $ \ti P_{\ti\mu} (z) = \( \frac1z - 1 \)^m =
(-1)^m z^{-m} (z-1)^m $.

Denote by $ \P_0 (\mu), \F_0 (\mu) $ subspaces $ \P_0, \F_0 $
appearing in the PaF structure $ \G_\mu = ( H, (\P_n), (\F_n), T )
$. That is, $ \F_0 \subset L_2 (\mu) $ is spanned by functions $ z
\mapsto z^k (1-z)^n $ for $ k \ge 0 $, and $ \P_0 (\mu) $ is spanned
by functions $ z \mapsto z^k \( \frac1z - 1 )^m $ or $ z \mapsto
z^{k-m} (z-1)^m $ for $ k \le 0 $.

The next lemma is well-known for finite measures (see
\cite[Sect.~XII.5, before Theorem 5.1]{Do}); here is a generalization
to moderate measures.

\begin{lemma}\label{2.2}
Let two functions, $ f $ on the circle and $ g $ on $ \R $, be related
by
\begin{equation}
f \bigg( \frac{\la-i}{\la+i} \bigg) = g (\la) \quad \text{for all }
\la \in \R \, .
\end{equation}
Then $ f \in \P_0 (\mu) $ if and only if $ g \in \P_0 (\nu) $. Also, $
f \in \F_0 (\mu) $ if and only if $ g \in \F_0 (\nu) $.
\end{lemma}

\begin{proof}
It suffices to prove the latter, $ f \in \F_0 (\mu) \equi g \in \F_0
(\nu) $, since $ f \in \P_0 (\mu) \equi \ti f \in \F_0 (\mu) $, and
$ g \in \P_0 (\nu) \equi \ti g \in \F_0 (\nu) $, and $ \ti f \(
\frac{\la-i}{\la+i} \) = f \( \frac{\la+i}{\la-i} \) = g(-\la) = \ti
g(\la) $.

In order to prove that $ f \in \F_0 (\mu) $ implies $ g \in \F_0 (\nu)
$, consider $ f(z) = z^k (1-z)^n $ for some $ k \ge 0 $; we have to
prove that $ g \in \F_0 (\nu) $, where $ g(\la) = \(
\frac{\la-i}{\la+i} \)^k \( \frac{2i}{\la+i} \)^n = (2i)^n \frac1{
(\la+i)^n } \( 1 - 2i \frac1{\la+i} \)^k $ is a linear combination of
functions $ \la \mapsto \frac1{ (\la+i)^{n+l} } $, $ l = 0,1,\dots,k
$. Such $ g $ is Fourier transform of a linear combination of
functions $ h_{n+l} (t) = t^{n+l-1} e^{-t} $ for $ t > 0 $ (otherwise
$ 0 $), except for the case $ n=l=0 $; in that case $ g $ is constant,
and we need Fourier transform of a measure (concentrated at the
origin) rather than a function of $ L_1 $. The same difficulty appears
for $ n > 0 $, when $ l = 0 $; in that case $ h_{n+l} = h_n $ does not
belong to $ W_1^n (\R) $, since $ h_n^{(n-1)} $ jumps at the origin,
and $ h_n^{(n)} $ is a finite measure rather than a function of $ L_1
$. However, a smoothing, say, $ t \mapsto \frac1\eps \int_{-\eps}^0
h_n (t+u) \, du $, does the job for $ l = 0 $. For $ l > 0 $ the
function $ h_{n+l} $ belongs to $ W_1^n (\R) $. So, $ g \in \F_0 (\nu)
$.

In order to prove that $ g \in \F_0 (\nu) $ implies $ f \in \F_0 (\mu)
$, consider $ g = \hat\phi $ where $ \phi \in W_1^n (\R) $, $ \phi (t)
= 0 $ for $ t \in (-\infty,0] $. The function $ \la \mapsto (\la+i)^n
g(\la) $ on the closed half-plane $ \Im \la \ge 0 $ is continuous, and
tends to $ 0 $ for $ |\la| \to \infty $. Therefore the function $ z
\mapsto (1-z)^{-n} f(z) $ on the closed disk $ |z| \le 1 $ is
continuous (and vanishes at $ 1 $). Take polynomials $ P_n $ such that
$ P_n (z) \to (1-z)^{-n} f(z) $ uniformly on the disk; then functions
$ z \mapsto (1-z)^n P_n (z) $ belong to $ \F_0 (\mu) $ and converge to
$ f $ in $ L_2 (\mu) $. So, $ f \in \F_0 (\mu) $.
\end{proof}

\section{Off-white noises}

Return to a Gaussian measure $ \ga $ in the space of distributions,
its spectral measure $ \nu $, and the corresponding stationary
Gaussian generalized random process. The spaces $ \P_0 (\nu), \F_0
(\nu) $ of $ L_2 (\nu) $, defined in Sect.~2, correspond unitarily
(via Fourier transform) to subspaces of the Hilbert space of all
$ \ga $-measurable linear functionals. Namely, $ \P_0 (\nu) $
corresponds to functionals localized (on the time axis) on $
(-\infty,0) $ (``the past''), and $ \F_0 (\nu) $ corresponds to
functionals localized on $ (0,\infty) $ (``the future''). Thus,
orthogonality of $ \P_0 (\nu), \F_0 (\nu) $ means independence of the
past and the future (which is the case for the white noise, whose
spectral measure is Lebesgue measure on $ \R $). The property
\begin{equation}\label{3.1}
L_2 (\nu) = \P_0 (\nu) \oplus \F_0 (\nu) \quad \text{in the FHS sense}
\end{equation}
means that dependence between the past and the future boils down to a
density. That is, $ \ga $ is equivalent (mutually absolutely
continuous) to another measure that makes the past and the future
independent.\footnote{%
 I mean $ (-\infty,0) $ and $ (0,\infty) $, not $ (-\infty,t) $ and $
 (t,\infty) $ for all $ t $ simultaneously.}
Such a process will be called an \emph{off-white noise.}

\begin{theorem}\label{3.2}
The following two conditions are equivalent.

(a) $ \nu $ is the spectral measure of an off-white noise;

(b) $ \mu $ has a density of the form $\displaystyle
\frac{d\mu}{d\mes} (z) = \bigg| \frac{ (z-z_1)\dots(z-z_{m-1}) }{
(z-1)^m } \bigg|^2 \exp \phi(z) $ for some $ m \in \{ 1,2,\dots \} $, some
points $ z_1,\dots,z_{m-1} $ on the circle, different from $ 1 $, and
some function $ \phi \in W_2^{1/2} $.
\end{theorem}

\begin{proof}
Condition (a) is equivalent to \eqref{3.1}. By Lemma \ref{2.2},
\eqref{3.1} is equivalent to $ L_2 (\mu) = \P_0 (\mu) \oplus \F_0
(\mu) $ in the FHS sense. The latter is \ref{prop2}(a) for $ N = -1
$. By Proposition \ref{prop2} it is equivalent to \ref{prop2}(b) for $
N = -1 $. It remains to note that $ \mu $ has no poles except for $ 1
$.
\end{proof}

Clearly, $ m $ in \ref{3.2}(b) is the same as $ m $ in \eqref{2.1}.
We are mostly interested in the case $ m = 1 $; in that case
\eqref{2.1} becomes
\begin{equation}\label{3.3}
\int_{-\infty}^{+\infty} \frac1{1+\la^2} \, \nu(d\la) < \infty \, ,
\end{equation}
and \ref{3.2}(b) becomes $ \frac{d\mu}{d\mes} (z) = \frac1{ |1-z|^2 }
w(z) $. However,\footnote{%
 I denote by $ \mes $ both Lebesgue measure on the circle and Lebesgue
 measure on $ \R $. We have $ dz = \frac{2i}{(\la+i)^2} \, d\la $,
 thus $ \mes (dz) = \frac2{\la^2+1} \mes(d\la) $; also, $ |1-z|^2 =
 \frac4{\la^2+1} $.}
$ \frac{d\nu}{d\mes} (\la) = \frac12 |1-z|^2
\frac{d\mu}{d\mes} (z) $ where $ z = \frac{\la-i}{\la+i} $; thus
\begin{equation}\label{3.4}
w \bigg( \frac{\la-i}{\la+i} \bigg) = 2 W (\la) \, ,
\end{equation}
where $ W (\la) = \frac{d\nu}{d\mes} (\la) $.

The condition $ \ln w \in W_2^{1/2} $ can be rewritten in terms of $ W
$,
\begin{equation}\label{3.6}
\iint \frac{ | \ln W(\la_1) - \ln W(\la_2) |^2 }{ | \la_1 - \la_2 |^2
} \, d\la_1 \, d\la_2 < \infty \, ,
\end{equation}
which is \eqref{eqSobolev} combined with the fact that $ \frac{dz_1
dz_2}{(z_1-z_2)^2} $ is invariant under linear-fractional
transformations.\footnote{%
 That is, if $ f(z) = \frac{az+b}{cz+d} $ then $ \frac1{(z_1-z_2)^2} =
 \frac{ f'(z_1) f'(z_2) }{ ( f(z_1) - f(z_2) )^2 } $.}
Recall that $ W(-\la) = W(\la) $.

\begin{proposition}
(a) If $ W $ satisfies \eqref{3.6} then
\begin{equation}\label{3.11}
\int_0^\infty | \ln W(2\la) - \ln W(\la) |^2 \frac{d\la}{\la} < \infty
\, .
\end{equation}
(b) Let $ W $ be strictly positive, have a continuous derivative, and
\begin{equation}\label{3.12}
\int_0^\infty \bigg| \frac{d}{d\la} \ln W(\la) \bigg|^2 \la \, d\la <
\infty \, .
\end{equation}
Then $ W $ satisfies \eqref{3.6}.
\end{proposition}

\begin{proof}
First, integration in \eqref{3.6} may be restricted from $ \R \times
\R $ to $ (0,\infty) \times (0,\infty) $. Indeed, using the property $
W(\la) = W(-\la) $ we get the kernel $ \frac2{|\la_1-\la_2|^2} +
\frac2{|\la_1+\la_2|^2} $ equivalent to $ \frac1{|\la_1-\la_2|^2} $.

Second,
\begin{multline}\label{3.13}
\int_0^\infty \int_0^\infty \frac{ | \ln W(\la_1) - \ln W(\la_2) |^2
 }{ | \la_1 - \la_2 |^2 } \, d\la_1 \, d\la_2 = \\
= \int_0^\infty
 \frac{du}{ (u-1)^2 } \int_0^\infty | \ln W(u\la) - \ln W(\la) |^2 \,
 \frac{d\la}{\la} \, ,
\end{multline}
which is just a change of variable, $ \la_1 = \la_2 u $.

\begin{sloppypar}
Let $ W $ satisfy \eqref{3.6}; we have to check \eqref{3.11}. Consider
$ f(u) = \bigg( \int_0^\infty | \ln W(u\la) - \ln W(\la)
|^2 \, \frac{d\la}{\la} \bigg)^{1/2} $. The triangle inequality gives
$ f(uv) \le f(u) + f(v) $, since $ \int_0^\infty | \ln W(uv\la) - \ln
W(v\la) |^2 \, \frac{d\la}{\la} = \int_0^\infty | \ln W(u\la) - \ln
W(\la) |^2 \, \frac{d\la}{\la} $. Also, $ f(u) < \infty $ for almost
all $ u $ due to \eqref{3.6} and \eqref{3.13}. Taking $ u $ such that
$ f(u) < \infty $ and $ f \( \frac2u \) < \infty $ we get $ f(2) <
\infty $, which is \eqref{3.11}.
\end{sloppypar}

\begin{sloppypar}
Let $ W $ satisfy \eqref{3.12}; we need to check \eqref{3.6}, or
equivalently, $ \int_0^\infty f^2(u) \frac{du}{(u-1)^2} < \infty $. We
have
\begin{gather*}
\bigg( \int_0^\infty \bigg| \int_1^u \frac{ W'(\la x) }{ W(\la x) } \,
 dx \bigg|^2 \la \, d\la \bigg)^{1/2} \le \int_1^u \bigg(
 \int_0^\infty \bigg| \frac{ W'(\la x) }{ W(\la x) } \bigg|^2 \la \,
 d\la \bigg)^{1/2} \, dx \, ; \\
\bigg( \int_0^\infty | \ln W(\la u) - \ln W(\la) |^2 \,
 \frac{d\la}{\la} \bigg)^{1/2} \le \bigg( \int_1^u \frac{dx}x \bigg) 
 \bigg( \int_0^\infty \bigg| \frac{ W'(\la) }{ W(\la) } \bigg|^2 \la
 \, d\la \bigg)^{1/2} \, ; \\
f(u) \le \const \cdot \ln u
\end{gather*}
for all $ u \in [1,\infty) $; similarly, $ f(u) \le \const \cdot |\ln
u| $ for all $ u \in (0,1] $. So, $ \int_0^\infty f^2(u) \,
\frac{du}{(u-1)^2} \le \const \cdot \int_0^\infty \( \frac{\ln u}{u-1}
\)^2 \, du < \infty $.
\end{sloppypar}
\end{proof}

\begin{example}
Assume that $ W(\la) = |\la|^\al $ for $ |\la| $
large enough, and $ W $ is strictly positive and smooth
everywhere. Then Condition \eqref{3.6} is satisfied if and only if $
\al = 0 $ (just the white noise).
\end{example}

\begin{example}
Assume that $ W(\la) = (\ln|\la|)^\al $ for
$ |\la| $ large enough, and $ W $ is strictly positive and smooth
everywhere. Then Condition \eqref{3.6} is satisfied for all $ \al $.
Condition \eqref{3.3} is also satisfied. Thus, every such
$ W $ describes an off-white noise.
\end{example}

\begin{example}
Assume that $ W(\la) = \exp ( -\ln^\al
|\la| ) $ for $ |\la| $ large enough (here $ \al > 0 $), and $ W $ is
strictly positive and smooth everywhere. Then Condition \eqref{3.6} is
satisfied if and only if $ \al < 1/2 $. Condition \eqref{3.3} is also
satisfied. So, for $ \al \in (0,1/2) $ every such $ W $ describes an
off-white noise.
\end{example}

\bigskip
\filbreak
{
\small
\begin{sc}
\parindent=0pt\baselineskip=12pt

School of Mathematics, Tel Aviv Univ., Tel Aviv
69978, Israel
\smallskip
\emailwww{tsirel@math.tau.ac.il}
{http://www.math.tau.ac.il/$\sim$tsirel/}
\end{sc}
}
\filbreak

\end{document}